\input amstex
\documentstyle{amsppt}
\magnification=1150

\def\ep{\varepsilon}

\nopagenumbers

\topmatter
\title A local regularity \\ of the complex Monge-Amp\`ere equation \endtitle
\rightheadtext{A local regularity of the complex Monge-Amp\`ere
equation}
\author Zbigniew B\l ocki and S\l awomir Dinew\endauthor
\leftheadtext{Zbigniew B\l ocki and S\l awomir Dinew}

\address \!\!Jagiellonian University, Institute of Mathematics,
\L ojasiewicza 6, 30-348 Krak\'ow, Poland \endaddress

\email Zbigniew.Blocki\@im.uj.edu.pl, Slawomir.Dinew\@im.uj.edu.pl
\endemail

\abstract We prove a local regularity (and a corresponding a
priori estmate) for plurisubharmonic solutions of the
nondegenerate complex Monge-Amp\`ere equation assuming that their
$W^{2,p}$-norm is under control for some $p>n(n-1)$. This
condition is optimal. We use in particular some methods developed
by Trudinger and an $L^q$-estimate for the complex Monge-Amp\`ere
equation due to Ko\l odziej.
\endabstract



\thanks Partially supported by the projects N N201 2683 35 and 189/6 PR
EU/2007/7 of the Polish Ministry of Science and Higher Education.
\endthanks

\endtopmatter

\document

\head Introduction \endhead

The aim of this note is to prove the following a priori estimate
for the complex Monge-Amp\`ere equation:

\proclaim{Theorem} Assume that $n\geq 2$ and $p>n(n-1)$. Let $u\in
W^{2,p}(\Omega)$, where $\Omega$ is a domain in $\Bbb C^n$, be a
plurisubharmonic solution of
  $$\det\big(u_{z_j\bar z_k}\big)=\psi>0.\tag 1$$
Assume that $\psi\in C^{1,1}(\Omega)$. Then for
$\Omega'\Subset\Omega$
  $$\sup_{\Omega'}\Delta u\leq C,$$
where $C$ is a constant depending only on $n$, $p$, $\text{\rm
dist}(\Omega',\partial\Omega)$, $\inf_\Omega\psi$,
$||\psi||_{C^{1,1}(\Omega)}$ and $||\Delta u||_{L^p(\Omega)}$.
\endproclaim

By a complex version of the Evans-Krylov theory (see e.g. [12] or
[5]), once one has an upper bound for the Laplacian (and thus for
mixed complex second derivatives) then also a
$C^{2,\alpha}$-estimate follows. We thus get the following local
regularity of plurisubharmonic solutions of (1)
  $$u\in W^{2,p}_{loc}\text{ for some }p>n(n-1),\ \psi\in C^\infty\
     \Longrightarrow\ u\in C^\infty.\tag 2$$
For $p>2n(n-1)$ this (and the theorem) is a consequence of a
general real theory from [14] (see [4]). For $p>n^2$ a similar a
priori estimate for $C^3$-solutions (without a regularity result
though) was recently shown in [8].

The main point about our result is that the condition $p>n(n-1)$
is essentially optimal. The fact that it is false for $p<n(n-1)$
follows from a complex counterpart of Pogorelov's example [11]
from [4]: the function
  $$u(z)=(1+|z_1|^2)|z'|^{2-2/n},$$
where $z'=(z_2,\dots,z_n)$, is in $W^{2,p}_{loc}$ for $p<n(n-1)$,
plurisubharmonic in $\Bbb C^n$, and satisfies
  $$\det\big(u_{z_j\bar z_k}\big)=c_n(1+|z_1|^2)^{n-2}\in
       C^\infty(\Bbb C^n)$$
($c_n$ is a constant depending only on $n$) in the weak sense of
[2].

The corresponding estimates and regularity for the real
Monge-Amp\`ere equation can be found in [15].

The main tool in the proof of Theorem 1 will be the following
estimate of Ko\l odziej [9] (see also [10]): if a plurisubharmonic
$u$ with $u\geq 0$ on $\partial\Omega$ solves (1) (with $\psi$
satisfying only $\psi\geq 0$) then for $q>1$ we have
  $$\sup_\Omega(-u)\leq C(q,n,\text{\rm diam}\,\Omega)\,
     ||\psi||_{L^q(\Omega)}^{1/n}.\tag 3$$
This result for $q=2$ is due to Cheng and Yau (see [1] and [7]).

\vskip 2mm

\noindent{\it Acknowledgement:} Part of the research was done
while the second named author was visiting the Princeton
University. He would like to thank this institution for the
perfect working conditions and hospitality and especially
professor Gang Tian for his encouragement and help.

\head Proof of Theorem 1 \endhead

By $C_1,C_2,\dots$ we will denote possibly different constants
depending only on the required quantities. Without loss of
generality we may assume that $\Omega=B$ is the unit ball in $\Bbb
C^n$ and that $u$ is defined in some neighborhood of $\bar B$. We
will use the notation $u_j=u_{z_j}$, $u_{\bar j}=u_{\bar z_j}$ and
$\Delta u=\sum_ju_{j\bar j}$. As usual, by $(u^{i\bar j})$ we will
denote the inverse transposed of $(u_{i\bar j})$.

We will first prove Theorem 1 assuming that $u$ is in $C^4$.
Differentiating (1) w.r.t. $z_p$ and $\bar z_p$ we will get
  $$u^{i\bar j}u_{i\bar jp}=(\log\psi)_p$$
and
  $$u^{i\bar j}u_{i\bar jp\bar p}=(\log\psi)_{p\bar p}
     +u^{i\bar l}u^{k\bar j}u_{k\bar l\bar p}u_{i\bar jp}.$$
Therefore
  $$u^{i\bar j}(\Delta u)_{i\bar j}\geq\Delta(\log\psi).
     \tag 4$$

We will now use an idea from [13]. For some $\alpha,\beta\geq 2$
to be determined later set
  $$w:=\eta(\Delta u)^\alpha,$$
where
  $$\eta(z):=(1-|z|^2)^\beta$$
Then
  $$w_i=\eta_i(\Delta u)^\alpha
     +\alpha\eta(\Delta u)^{\alpha-1}(\Delta u)_i$$
and
  $$\align u^{i\bar j}w_{i\bar j}
     &=\alpha\eta(\Delta u)^{\alpha-1}u^{i\bar j}(\Delta u)_{i\bar j}
        +\alpha(\alpha-1)\eta(\Delta u)^{\alpha-2}
              u^{i\bar j}(\Delta u)_i(\Delta u)_{\bar j}\\
     &\ \ \ \ \ \ \ \ +2\alpha(\Delta u)^{\alpha-1}\text{Re}
               \big(u^{i\bar j}\eta_i(\Delta u)_{\bar j}\big)
         +(\Delta u)^\alpha u^{i\bar j}\eta_{i\bar j}.\endalign$$
By (4) and the Schwarz inequality for $t>0$
  $$\align u^{i\bar j}w_{i\bar j}
     &\geq\alpha\eta(\Delta u)^{\alpha-1}\Delta(\log\psi)
        +\alpha(\alpha-1)\eta(\Delta u)^{\alpha-2}
              u^{i\bar j}(\Delta u)_i(\Delta u)_{\bar j}\\
     &\ \ \ \ \ \ \ \ -t\alpha(\Delta u)^{\alpha-1}
          u^{i\bar j}(\Delta u)_i(\Delta u)_{\bar j}
        -\frac 1t\alpha(\Delta u)^{\alpha-1}u^{i\bar j}\eta_i\eta_{\bar j}
         +(\Delta u)^\alpha u^{i\bar j}\eta_{i\bar j}.\endalign$$
Therefore with $t=(\alpha-1)\eta/\Delta u$ we get
  $$u^{i\bar j}w_{i\bar j}
      \geq\alpha\eta(\Delta u)^{\alpha-1}\Delta(\log\psi)
         +(\Delta u)^\alpha u^{i\bar j}\left(\eta_{i\bar j}
            -\frac\alpha{\alpha-1}\frac{\eta_i\eta_{\bar j}}\eta\right).$$
We now have
  $$\align\eta_i&=-\beta z_i\eta^{1-1/\beta}\\
      \eta_{i\bar j}&=-\beta\delta_{i\bar j}\eta^{1-1/\beta}
         +\beta(\beta-1)\bar z_iz_j\eta^{1-2/\beta}\endalign$$
and thus
  $$\big|\eta_{i\bar j}\big|,\ \big|\frac{\eta_i\eta_{\bar j}}\eta\big|\leq
      C(\beta)\,\eta^{1-2/\beta}.$$
We will get
  $$u^{i\bar j}w_{i\bar j}\geq -C_1(\Delta u)^{\alpha-1}
         -C_2w^{1-2/\beta}(\Delta u)^{2\alpha/\beta}\sum_{i,j}|u^{i\bar j}|.$$

Fix $q$ with $1<q<p/(n(n-1))$. Since $||\Delta u||_p$ (this way we
will denote norms in $L^p(B)$) is under control, it follows that
$||u_{i\bar j}||_p$ and $||u^{i\bar j}||_{p/(n-1)}$ are as well.
It follows that for
  $$\alpha=1+\frac p{qn},\ \ \ \beta=2\big(1+\frac{qn}p\big)$$
we have
  $$||(u^{i\bar j}w_{i\bar j})_-||_{qn}
       \leq C_3(1+(\sup_B\,w)^{1-2/\beta}),$$
where $f_-:=-\min(f,0)$.

By [2] we can find continuous plurisubharmonic $v$ vanishing on
$\partial B$ and such that
  $$\det(v_{i\bar j})=((u^{i\bar j}w_{i\bar j})_-)^n$$
(weakly). Essentially by an inequality between arithmetic and
geometric means (see [3] how to extend it to the weak case) we
have
  $$\align u^{i\bar j}v_{i\bar j}
     &\geq n(\det(u^{i\bar j}))^{1/n}(\det(v_{i\bar j}))^{1/n}\\
     &=n\psi^{-1/n}(u^{i\bar j}w_{i\bar j})_-\\
     &\geq-\frac 1{C_4}u^{i\bar j}w_{i\bar j}.\endalign$$
It follows that $w\leq -C_4v$ and by Ko\l odziej's inequality (3)
  $$\align\sup_Bw&\leq C_5||\det(v_{i\bar j})||_q^{1/n}\\
     &=C_5||(u^{i\bar j}w_{i\bar j})_-||_{qn}\\
     &\leq C_6(1+(\sup_B\,w)^{1-2/\beta}).\endalign$$
Therefore $w\leq C_7$ and the desired estimate follows if $u\in
C^4$.

Now assume that the solution is just in $W^{2,p}$. Similarly to
[2], instead of $\Delta u$ we will consider for $\ep>0$ the
following approximations to the Laplacian
  $$T=T_\ep u=\frac{n+1}{\ep^2}(u_\ep-u),$$
where
  $$u_\ep(z)=\frac 1{\lambda(B(z,\ep))}\int_{B(z,\ep)}u\,d\lambda$$
and $\lambda$ denotes the Lebesgue measure in $\Bbb C^n$. Since
$T_\ep u\to\Delta u$ weakly as $\ep\to 0$, it is enough to show a
uniform upper bound for $T$ independent of $\ep$.

By [2] we have
  $$u^{i\bar j}u_{\ep,i\bar j}
     \geq n\psi^{-1/n}(\det(u_{\ep,i\bar j}))^{1/n}
     \geq n\psi^{-1/n}(\psi^{1/n})_\ep$$
and thus, coupling this with $u^{i\bar j}u_{i\bar j}=n$, we obtain
the following counterpart of (4)
  $$u^{i\bar j}T_{i\bar j}\geq n\psi^{-1/n}T_\ep(\psi^{1/n})
     \geq -C_8.$$
Changing the definition of $w$ to $\eta T^\alpha$ (since $u$ is
plurisubharmonic, $T$ is nonnegative, hence $T^{\alpha}$ is well
defined)  and repeating the previous computations we will get
  $$u^{i\bar j}w_{i\bar j}\geq C_9T^{\alpha-1}
         -C_{10}w^{1-2/\beta}T^{2\alpha/\beta}\sum_{i,j}|u^{i\bar j}|.$$
The rest of the proof is now the same as before.

\enddocument